\documentclass[pdflatex,sn-mathphys-num]{sn-jnl}

\usepackage{mathrsfs}

\usepackage{graphicx}%
\usepackage{multirow}%
\usepackage{amsmath,amssymb,amsfonts}%
\usepackage{amsthm}%
\usepackage{mathrsfs}%
\usepackage[title]{appendix}%
\usepackage{xcolor}%
\usepackage{textcomp}%
\usepackage{manyfoot}%
\usepackage{booktabs}%
\usepackage{algorithm}%
\usepackage{algorithmicx}%
\usepackage{algpseudocode}%
\usepackage{listings}%

\theoremstyle{thmstyleone}%
\newtheorem{theorem}{Theorem}
\theoremstyle{thmstyletwo}%

\theoremstyle{thmstylethree}%
\newtheorem{definition}{Definition}%

\begin{document}

\title[Article Title]{A distribution-guided Mapper algorithm}


\author[1]{\fnm{Yuyang} \sur{ Tao}}\email{taoyy2022@shanghaitech.edu.cn}

\author*[1]{\fnm{Shufei} \sur{ Ge}}\email{geshf@shanghaitech.edu.cn}

\affil*[1]{\orgdiv{Institute of Mathematical Sciences}, \orgname{ShanghaiTech University}, \orgaddress{\street{393 Middle Huaxia Road}, \postcode{201210}, \state{Shanghai}, \country{Country}}}




\abstract{
\textbf{Background:} The Mapper algorithm is an essential tool for exploring the data shape in topological data analysis. With a dataset as an input, the Mapper algorithm outputs a graph representing the topological features of the whole dataset. This graph is often regarded as an approximation of a Reeb graph of a dataset. The classic Mapper algorithm uses fixed interval lengths and overlapping ratios, which might fail to reveal subtle features of a dataset, especially when the underlying structure is complex.\\
\textbf{Results:} In this work, we introduce a distribution-guided Mapper algorithm named D-Mapper, which utilizes the property of the probability model and data intrinsic characteristics to generate density-guided covers and provide enhanced topological features. Moreover, we introduce a metric accounting for both the quality of overlap clustering and extended persistent homology to measure the performance of Mapper-type algorithms.  Our numerical experiments indicate that the D-Mapper outperforms the classic Mapper algorithm in various scenarios. We also apply the D-Mapper to a SARS-COV-2 coronavirus RNA sequence dataset to explore the topological structure of different virus variants. The results indicate that the D-Mapper algorithm can reveal both the vertical and horizontal evolutionary processes of the viruses. Our code is available at \url{https://github.com/ShufeiGe/D-Mapper}.\\
\textbf{Conclusion:} The D-Mapper algorithm can generate covers from data based on a probability model. This work demonstrates the power of fusing probabilistic models with Mapper algorithms.\\}





\maketitle

\section{Introduction}\label{sec1}

In recent decades, machine learning methods have become popular tools to discover valuable information from data. These methods typically focus on prediction tasks by establishing a mapping between responses and predictors. The data shape may reveal important information in data from a new perspective and provide alternative insights to supplement prediction tasks, which is often overlooked in the learning of mapping. To be specific, the data shape refers to the manifold formed by the support set of data distribution. The data shape reveals the data distribution, reflecting the spatial correlations and dependency structures among data points. Such information is crucial in clustering and feature extraction tasks. However, many existing methods cannot fully utilize such information.

Topology data analysis (TDA) is a specific area of investigating data shapes based on topological theory. Topology is a powerful tool for studying the shape of objects by finding invariants that remain unchanged under continuous deformations. These invariants can reflect objects' intrinsic features. TDA attempts to identify and utilize these topology invariants in data analysis. The Mapper algorithm is a simple but powerful TDA tool for visualizing and clustering data. The algorithm constructs a graph in which each vertex represents a cluster, and each edge represents the two adjacent clusters that share some elements. The colour and size of vertices are defined by the users and can provide additional information. For example, the size of vertices and the colour depth can represent the average value and the number of elements in the cluster. 
The Mapper algorithm was first proposed in \cite{singhTopologicalMethodsAnalysis} and has been applied in various domains, such as social media \cite{almgrenMiningSocialMedia2017}, fraud detection \cite{mitraExperimentsFraudDetection2021}, and evolutionary computation \cite{todaVisualizationClusteringGraph2022}. The Mapper algorithm is especially suitable for biology data, which are usually complex and high-dimensional  \cite{skafTopologicalDataAnalysis2022}. For example, the Mapper algorithm was applied to breast cancer transcriptional microarray data and successfully identified two subgroups with $100\%$ survival rate \cite{nicolauTopologyBasedData2011}. It is also used to analyze transcriptional programs that control cellular lineage commitment and differentiation during development, and the proposed scTDA identified four transient states over time \cite{rizviSinglecellTopologicalRNAseq2017}. Exploring the conformation space of proteins is another important task in computational biology, and the Mapper algorithm has been successfully applied to analyze intermediate conformations, the results are closely consistent with experimental findings \cite{dehghanpoorUsingTopologicalData}. 

The classic Mapper algorithm requires a filter function $f: X \rightarrow \mathbb{R}^{p}$ to project data $X$ onto the Euclidean space $\mathbb{R}^{p}$ and a set of open cover $\mathcal{U}$ on the projected data. In most of the literature, $p$ is set to $1$, i.e.  the filter function $f$ maps data from a high dimensional space to a real line. Different filter functions and open covers may result in different outputs, so it is necessary to select them carefully. Improper filter functions or open covers may not accurately reveal data shapes, resulting in poor overlap  clustering. A filter function is chosen primarily based on the features of interest and also depends on specific applications. 
 Determining the optimal parameters of the model, such as the overlapping rates and interval lengths, usually requires extensive manual tuning and experiments. As pointed out in \cite{carriereStatisticalAnalysisParameter2017}, cover choices lead to a broad range of possible Mapper outputs, resulting in many different graph depictions and clustering results. Unreasonable choices may cause the disappearance of certain topological features.  

Many works have been proposed to improve the performance of the classic Mapper. For example, \cite{dlotkoBallMapperShape2019} proposed to generate covers for the dataset by constructing a set of balls directly. This method saved the trouble of choosing the filter function.  \cite{buiFMapperFuzzyMapper2020} used a fuzzy clustering algorithm to generate covers automatically with random overlapping ratios.  To achieve adaptive cover construction, an information criteria was developed based on an X-means algorithm to generate adaptive covers \cite{chalapathiAdaptiveCoversMapper2021}.

In this work, alternatively, we propose a distribution-guided Mapper (D-Mapper) algorithm to relax the restriction of regular covers and overlapping ratios. Our proposed algorithm utilizes the property of the probability model and data intrinsic characteristics to generate distribution-guided covers and provides enhanced topological features. 
 Unlike the classic Mapper algorithm, our proposed D-Mapper does not rely on predefined static overlapping rates and interval lengths. Instead, we fit the projected data to a mixture distribution model and automatically generate flexible covers. Moreover, evaluating the performance of the Mapper algorithm is challenging due to the complexity of the output graph structure. Most methods in literature only evaluate it from a clustering perspective, overlooking the quality of topological structure~\cite{buiFMapperFuzzyMapper2020, chalapathiAdaptiveCoversMapper2021}. In this work, we introduce a metric that can quantitatively and objectively evaluate the performance of the Mapper algorithm from both clustering and topological aspects. We use the extended persistent homology as a tool to capture topological signatures, and distinguish noise from signal on the diagram by constructing a confidence set on the diagram. The signal rate can be used as a metric to evaluate the quality of topological structure. The simulation studies indicate that the D-Mapper algorithm outperforms the classic Mapper algorithm.  The remainder of the paper is organized as follows. We briefly review the classic Mapper algorithm and introduce our D-Mapper and evaluation metric in the section \ref{Methods}. Section \ref{Results} focuses on comparing the classic Mapper with the D-Mapper on simulated datasets and applying the proposed D-Mapper algorithm to a real-world SARS-COV-2 RNA sequence dataset. We discuss and summarize our work in Sections \ref{Discussion} and \ref{Conclusion}.

\section{Methods}\label{Methods}

\subsection{Basic notions}\label{subsec2.1}
Before we introduce our method, we first briefly review the essential background of the Mapper algorithm. The theoretical foundation of the Mapper algorithm is the nerve theorem which guarantees the nerves produced by a cover on a topological space $X$ is homotopy equivalent to that space. We introduce the following notions to describe the nerve  theorem \cite{carlssonTopologyData2009,chazalIntroductionTopologicalData2021}.

\begin{definition}[Simplex]
Given a set $P = \{p_0,...,p_k\} \subset \mathbb{R}^{p} $ of $k+1$ affinely independent points, a $k$-dimensional simplex $\sigma$, or $k$-simplex for short, spanned by $P$ is the set of convex combinations such that:
$$\{ x| x= \sum_{i=0}^k \lambda_i p_i, \sum_{i=0}^k \lambda_i = 1 ,\lambda_i \geq 0   \}.$$
And the points of $P$ are the vertices of simplex $\sigma$ and the simplices spanned by the subsets of P are the faces of simplex $\sigma$.
\end{definition}

\begin{definition}[Geometric simplicial complex]
A geometric simplicial complex $K$ in $\mathbb{R}^{n}$ is a finite collection of simplices satisfying the following two conditions:\\
a) Arbitrary face of any simplex of $K$ is a simplex of $K$.\\
b) The intersection of any two simplices of $K$ is either empty or a common face of the two.
\end{definition}
The union of the simplices of $K$ constitutes the underlying space of $K$, denoted as $|K|$, which inherits from the topology space of $\mathbb{R}^{p}$. Thus, the geometric simplicial complex can also be regarded as a topological space.

\begin{definition}[Abstract simplicial complex]
Let $V$ be a finite set. An abstract simplicial complex $\mathcal{K}$ given the set $V$ is a set of finite subsets of $V$ such that:\\
a) All elements of  $V$ belongs to $\mathcal{K}$.\\
b) If $\tau \in \mathcal{K}$ , any subset of $\tau$ belongs to $\mathcal{K}$.
\end{definition}
 
\begin{definition}[Open cover]
Suppose  $\mathcal{U} = (u_i) ,i \in I$ is a collection of open subset of a topological space $X$, then we say  $\mathcal{U}$ is an open cover of $X$ if black $X = \bigcup_{i \in I} u_i $. 
\end{definition}
Given an open cover of topological space $X$, $\mathcal{U} = (u_i) $, the nerve of $\mathcal{U}$ is an abstract simplicial complex C$(\mathcal{U})$ with vertex set $\mathcal{U}$. With these definitions, we can introduce the most important theorem in constructing the Mapper algorithm.
\begin{theorem}[Nerve theorem]
 Let $\mathcal{U} = (u_i),~i \in I$ be an open cover of a  paracompact space $X$ by open sets such that the intersection of any sub-collection of the $u_i$’s is either empty or contractible. Then, $X$ and the nerve C($\mathcal{U}$) are homotopy equivalent.
\end{theorem}

This theorem allows us to map the topology of continuous into abstract combinatorial structures by building a nerve complex. It bridges the gap between continuous space and its discrete representation.  This crucial theorem is also a motivation for the construction of the Mapper algorithm.

\subsection{Mapper algorithm}
The Mapper algorithm is an important tool to construct the discrete version of Reeb graphs, which encode connected information of the support manifold \cite{dey2017topological}. Algorithm \ref{Mapper} depicts the classic Mapper algorithm. Firstly, the original data is projected onto a real line by a user-specified filter function (Algorithm \ref{Mapper}, line 1). To construct a cover $\mathcal{U}$ on the projected data, the number of components $n$ and overlapping ratio $p$ should be chosen carefully (Algorithm \ref{Mapper}, line 2). The improper choice of these parameters may lead to failure estimation of the data shape. With these two parameters, a cover with $n$ regular intervals can be obtained (Algorithm \ref{Mapper}, line 3). Then, the inverse image of the cover $\mathcal{U}$ is achieved on the original data, $f^{-1}(\mathcal{U})$, generating some hypercubes in the original data space. This process is often called pulling back (Algorithm \ref{Mapper}, line 4). Finally, clustering data points within each hypercube. Each cluster corresponds to a vertex in the Mapper graph. If two vertices share any elements, an edge is added between these two  vertices (Algorithm \ref{Mapper}, line 5).
\begin{algorithm}[H]
    \caption{The classic Mapper algorithm}
    \label{Mapper}
    \begin{algorithmic}[1]
      \State Choose a proper filter function $f$ to project data onto a real line, $f: X \rightarrow \mathbb{R}$.
      \State Choose a component number $n$ and overlapping percentage ratio $p$.
      \State Construct a cover $\mathcal{U} = (u_i), i=1...n$ on the projected data $f(X)$ based on the parameters $n$ and $p$.
      \State Pull back the intervals of the projected data to the original space, $f^{-1}(\mathcal{U})$. 
      \State Cluster on the refined cover and build the nerves with the clustering result.
    \end{algorithmic}
\end{algorithm}

The classic Mapper algorithm is simple but powerful for exploring and visualizing data, while the selection of parameters $n, ~p$ involves extensive manual tuning, and the algorithm is sensitive to these parameters. In addition, the Mapper algorithm's flexibility is often restricted by regular intervals and fixed overlapping ratios, which may hinder the discovery of complex data structures.

\subsection{D-Mapper}
The restriction of regular intervals with fixed overlapping ratios is one of the major limitations of the Mapper algorithm. We propose a distribution-guided Mapper algorithm to generate flexible covers to reflect the underlying data structures.   Our proposed method automatically chooses the overlapping ratios based on the distribution of the projected data and produces more flexible covers to reveal the data shapes more accurately. The key idea of our algorithm is to fit the projected data with a mixture probability model.  Each component in the mixture model can be viewed as an interval, and the probability (likelihood) of each data point assigned to each interval can be explicitly calculated.  Once we get the mixture distribution of the projected data, we can create intervals based on the distribution in many ways. 
 
\begin{figure}[!t]%
\centering
    \includegraphics[width=0.9\textwidth]{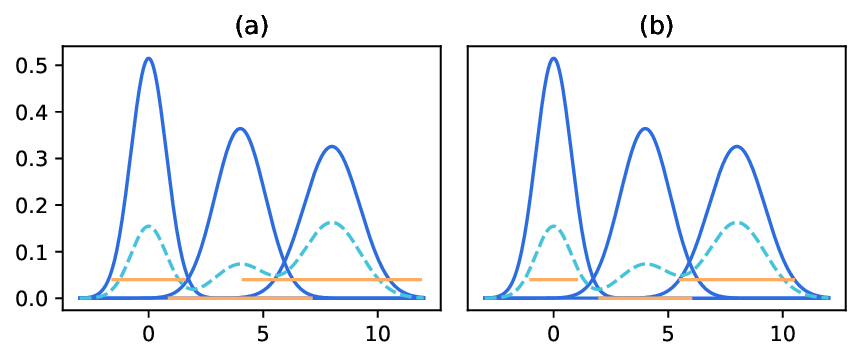}
    \vspace{-4mm}
    \caption{An illustration of intervals produced by the D-Mapper algorithm. The deep blue line represents the probability density function of each component in the GMM. The shallow blue dashed line presents the probability density function of the GMM. The orange lines are intervals that are produced naturally given a confidence level of $1-\alpha$. (a) When $\alpha = 0.01$, there are overlaps between adjacent intervals. (b) When $\alpha = 0.1$, there is a gap between the first and second intervals.  $\alpha$ controls the overlap of intervals, and it should be chosen carefully.} 
    \label{show_algorithm_alpha}
\end{figure}

 Here, we introduce a simple way to construct intervals based on a confidence level of $1-\alpha$ for a given probability distribution.  With a proper selection, $1-\alpha$ confidence intervals of each component of the mixture model can automatically produce some overlaps. Figure \ref{show_algorithm_alpha} shows the idea of naturally producing intervals given a confidence level of $1-\alpha$. This attribute provides a natural scheme for constructing flexible covers on the projected data. The specific procedures are as follows.\\
1) Choose an appropriate number of intervals $n$ and a confidence level of  $1-\alpha$. The number of components naturally matches the number of intervals. \\
2) Use a mixture model to fit the projected data.\\
3) Intervals are determined by $1-\alpha$ confidence intervals of each component of the mixture model.

In both our simulation and real data experiments, we implement the Gaussian mixture models (GMM) to fit the projected data due to its simplicity and flexibility. The model inference is done via the expectation maximization (EM) algorithm \cite{2009The}. Other mixture models or any distribution with multiple modes can be used as alternatives to the GMM. 

By incorporating a mixture model into the Mapper algorithm, it can produce flexible intervals and overlaps. We call this distribution-guided Mapper algorithm the D-Mapper algorithm. With these flexible intervals, we can construct the output graph similar to the classic Mapper algorithm: pulling back these intervals to the original data space, then clustering data points within each interval separately, all sub-groups constitute the vertices of the nerve, and adding an edge between two vertices if there are any shared data points. Algorithm \ref{D-Mapper} gives details of the D-Mapper.

\begin{algorithm}[H]
    \caption{The D-Mapper aglorithm}
    \label{D-Mapper}
    \begin{algorithmic}[1]
      \State Choose a proper filter function $f$ to project data onto a real line,  {$f: X \rightarrow \mathbb{R}$}. 
      \State Choose a component number $n$ and a confidence level $1-\alpha$.
      \State Fit the projected data to a mixture model.
      \For{$i$th component of the mixture model} 
        \State Set the $\frac{\alpha}{2}$ quantile $s_i$ as the start point of the interval,
        \State Set the $1-\frac{\alpha}{2}$ quantile $e_i$ as the end point of the interval,
        \State The interval of $i$th component is $u_i = (s_i,e_i)$.
      \EndFor
      \State The collection of intervals $\mathcal{U} = (u_i) ,i=1...n,$ is a cover on the projected data $f(X)$.
      \State Pull back the intervals of the projected data to the original space, $f^{-1}(\mathcal{U})$.
      \State Cluster on the refined cover and build the nerve with the clustering result.
    \end{algorithmic}
\end{algorithm}

The confidence level $1-\alpha$ controls the overlapping ratios:  a larger $\alpha$ value leads to lower overlapping ratios.  The distributions of components in mixture models are often heterogeneous, automatically resulting in flexible intervals. Proper $\alpha$  selections are crucial to ensure all points are covered.  The Mapper algorithms usually require pairwise overlapping (i.e. each interval overlaps with its neighbours)\cite{dey2017topological,dlotkoBallMapperShape2019,buiFMapperFuzzyMapper2020,chalapathiAdaptiveCoversMapper2021}. The D-Mapper algorithm can also reserve the pairwise overlap property, as shown in Figure \ref{show_algorithm_alpha} (a).  However, a larger $\alpha$ value may result in disjoint intervals, Figure \ref{show_algorithm_alpha} (b) gives an example of disjoint intervals caused by an improper large $\alpha$ value. Thus, in this work, we propose a method to find the upper bound of $\alpha$ that guarantees the pairwise overlap property. By establishing this upper bound, we also significantly narrow the potential range of $\alpha$ values, thereby facilitating parameter tuning.  Denote $F_i^{-1}$ the inverse of the cumulative density function of the $i$th ordered component, thus the $i$th interval $u_i$ is given by $[F_i^{-1}(\alpha/2),F_i^{-1}(1-\alpha/2)]$, $i=1,\ldots,n$.  The idea of this method is to find the upper bound $\alpha'$, that guarantees the intersection of each paired adjacent intervals is not empty, $u_i\cap u_{i+1}\neq \emptyset$, \emph{i.e.}, $F_i^{-1}(1-\alpha/2) \ge F_{i+1}^{-1}(\alpha/2)$, $i=1,2,\ldots,n-1$. Notice that although pairwise overlapping is a good property for a Mapper graph, this property is not necessary for all situations. One can construct a good Mapper graph with some paired adjacent intervals being non-overlap as long as the union of intervals can cover all data points. As shown in Figure \ref{cirs_results} (a) in Section \ref{Results}, since there are no points between the two disjoint circles, disjoint paired adjacent intervals should be presented in the corresponding Mapper output graph. We add one threshold parameter $\alpha^*$ to allow for disjoint paired adjacent intervals (lines 6-10 of Algorithm \ref{alpha-upperbound}). In practice, we suggest to set $\alpha^{*}=0.005$.  Algorithm \ref{alpha-upperbound} describes how to find the upper bound of $\alpha$. 

\begin{algorithm}[H]
    \caption{$\alpha$ upper bound}
    \label{alpha-upperbound}
    \begin{algorithmic}[1]
      \State $\boldsymbol{J}$ =  $\emptyset$
      \For{$i$ in range  (0, $n-1$)} 
        \State $S_1(\alpha) = F_i^{-1}(1-\alpha/2)$,
        \State $S_2(\alpha) =  F_{i+1}^{-1}(\alpha/2)$,
        \State Let $S_1(\alpha) = S_2(\alpha)$, solve for solution $\alpha_i'$.
      \If {$\alpha_i' \geq \alpha^*$}
        \State $ \boldsymbol{J} = \boldsymbol{J}\cup \{i\}$.
      \EndIf
      \EndFor
      \State $\alpha' = min\{\alpha_i'\}_{i \in \boldsymbol{J}}$.
      \State The upper bound of $\alpha$ is $\alpha'$, $\alpha \in (0,\alpha')$.
    \end{algorithmic}
\end{algorithm}

\subsection{Evaluation metric}
The silhouette coefficient ($SC$) is often used as a measure to evaluate the quality of (overlap) clustering \cite{hanClusterAnalysis2012}. The $SC$ assesses how well the clusters are separated, how close the clusters are, and whether they are stable for overlap clustering. Suppose we have $n$ samples in dataset $D$ that could be divided into $k$ clusters: $C_1,...,C_k$. For any two data points $x,~x'\in D$, the compactness of the cluster to which $x$ belongs can be defined as
\begin{equation}
a(x) = \frac{\sum_{x' \in C_i, x \neq x'} d(x,x')}{|C_i|-1},
\end{equation}

where $d(x,x')$ represents the distance between $x$ and $x'$, $|C_i|$ is the number of data points in cluster $i$. The degree of  separation between $x$ and other clusters can be computed as
\begin{equation}
    b(x) = \mathop{\min}_{C_j:j \neq i} \bigg \{ \frac{\sum_{x' \in C_j} d(x,x')}{|C_j|}  \bigg \}.
\end{equation}
And the $SC$ of $x$ is

\begin{equation}
SC(x) = \frac{b(x)-a(x)}{\max \{a(x),b(x)\}}.
\end{equation}

The value of the $SC$ ranges from $-1$ to $1$. A value close to $1$ indicates the point is close to the current cluster, and more distinct from other clusters. If the $SC$ is less than $0$, it means the point is closer to other clusters compared to the current cluster, and it usually indicates bad clustering results. The clustering of the whole dataset can be assessed by taking the average of the $SC$s of all points.

\begin{figure}[t]%
\centering
    \includegraphics[width=0.9\textwidth]{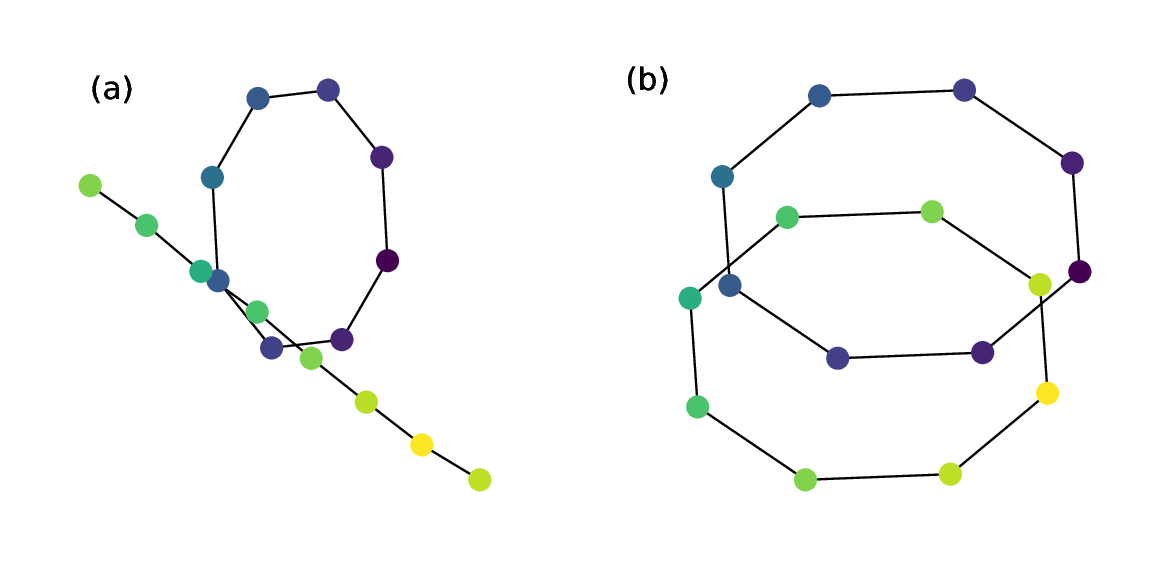}
    \caption{An example of the classic Mapper algorithm on the same dataset but outputs different graphs. The dataset is shown in Figure \ref{cirs_results} (a), the classic Mapper is implemented,  and the clustering algorithm is DBSCAN with a radius of $0.5$ and a minimum of samples $3$. (a) The output graph when $n=12, p=0.01$ and its $SC = 0.283$, $SC_{adj}=0.521$. This graph has a higher $SC$ but a poor topological structure. (b) The output graph when $n=12, p=0.1$ and its $SC = 0.246$, $SC_{adj}=0.812$. This graph has a lower $SC$ but a good topological structure. }
    \label{show_case}
\end{figure}

The $SC$ reflects the clustering quality of data points, but it does not evaluate the topology structure of a Mapper output graph. Figure \ref{show_case} shows an example of different Mapper graphs on the same dataset. Figure \ref{show_case} (a)  has a higher $SC$ but a poor topological feature. Conversely, Figure \ref{show_case} (b)  has a lower $SC$ but a good topological structure.  Figure \ref{show_case} (a) shows that a small overlap $p$ makes one of the circles disconnected, and however, more separable clusters lead to a high  $SC$. Therefore, a metric that evaluates both clustering and topological structure is needed. Although, the topological information encoded in  Mapper has been well studied theoretically \cite{dey2017topological}. There is no practical metric to evaluate the topology structure of Mapper graphs. Most existing work evaluates the topology structure of Mapper graphs based on the clustering results only. The extended persistence diagram has been proven to be a powerful tool for capturing topological signatures of Mapper graph \cite{carriereStructureStabilityOneDimensional2018}. In this manuscript, we alleviate it to provide a simple quantitative metric for evaluating the topology structure of Mapper output.

{The extended persistence diagram is an extension of the persistence diagram, the extended persistence diagram uses both sub-level sets and super-level sets as a filtration, therefore it contains more information than the original persistence diagram \cite{Cohen-Steiner2009}. In practice}, given a graph with a function defined on its nodes, we can compute the extended persistence diagram that reflects the topological features of this graph. For a Mapper graph, the function on a node can be naturally defined as the mean value of points in the node. The points on the diagram can be regarded as signatures of the Mapper graph. However, these signatures are not always meaningful, points near the diagonal are often seen as noise. To distinguish the real signal from the noise, one simple yet effective method is to calculate the confidence interval using bottleneck bootstrap and then identify the noise based on this interval \cite{carriereStatisticalAnalysisParameter2017}. The bottleneck bootstrap is an effective way to compute confidence set on a persistence diagram \cite{fasyConfidenceSetsPersistence2014}. It uses bootstrap samples to get a bootstrap Mapper graph and then calculates the bottleneck distance from the original Mapper graph. Repeat this step many times, and an approximate distribution of bottleneck distances can be obtained, then we can easily calculate the confidence set with this approximate distribution. See Algorithm \ref{confidence_sets} in Appendix \ref{bootstrap} for more details of the bottleneck bootstrap algorithm.

We introduce a metric to evaluate the quality of the extended persistence diagram depending on the confidence set; we call this metric the topological signal rate, denoted as $TSR$. The $TSR$ is defined as the number of real signal points divided by the total number of points on the extended persistence diagram, serving as a quantitative indicator of the quality of the extended persistence diagram. 

\begin{definition}[Topological signal rate]
\label{tsrdef}
The topological signal rate is a scalar, it evaluates the quality of a persistence diagram or extended persistence diagram, denoted as $TSR$,
$$TSR = \frac{N_{signal}}{N},$$
where $N$ is the number of all points on the diagram and $N_{signal}$ is the number of topological signals on the diagram. 
\end{definition}

\begin{figure}[h]%
\centering
\includegraphics[width=0.7\textwidth]{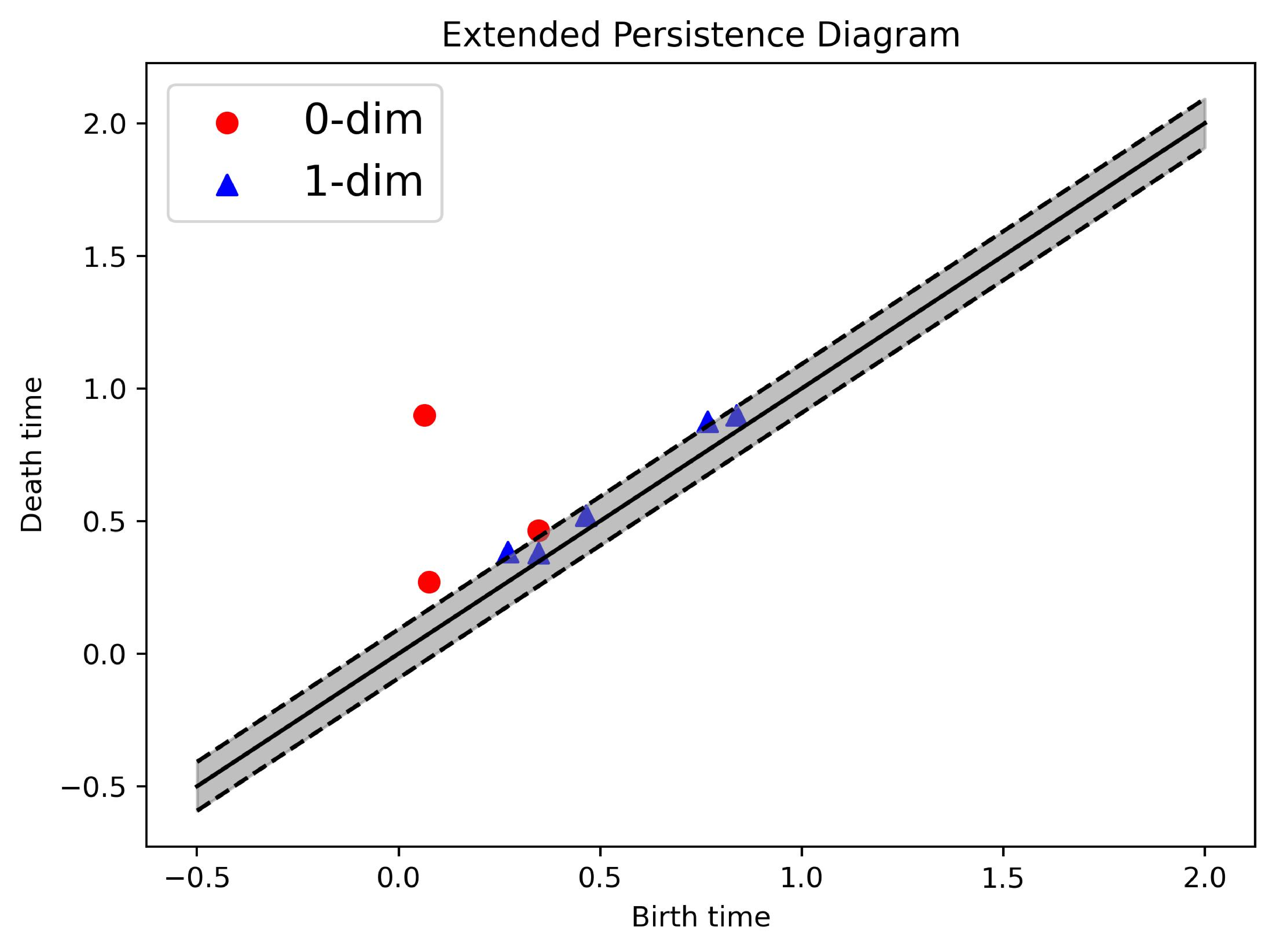}
\vspace{-2mm}
\caption{An example of the extended persistence diagram. There are a total of $8$ points in the diagram; the gray area is computed by the bottleneck bootstrap, and points inside this area are noise. Thus, $6$ points are noises, $2$ points are signals, and the $TSR$ is $0.25$.}
\label{EX_PH}
\end{figure}

An example of an extended persistence diagram and its $TSR$ is illustrated in Figure \ref{EX_PH}. The gray area is the confidence set estimated by the bottleneck bootstrap. Points inside are noises, and those outside are signal points. The $TSR$ in this example is $0.25$, indicating poor quality of the corresponding Mapper graph. To get a proper evaluation of the Mapper graph, we can combine the $TSR$ with $SC$ through weighted averaging. We call this metric an adjusted silhouette coefficient.

\begin{definition}[Adjusted silhouette coefficient]
The adjusted silhouette coefficient is the weighted average of the normalized $SC$ and the $TSR$, denoted as $SC_{adj}$, this value can serve as a metric to evaluate the Mapper-type algorithms:
$$SC_{adj} = w_1 SC_{norm}+w_2 TSR,$$
where $TSR$ is given by Definition \ref{tsrdef}, representing the topological signal rate, $SC_{norm} = \frac{SC+1}{2}$ is the normalized value of $SC$, and we choose $w_1 = w_2 = 0.5$ here to give same weights of clustering and topological structure.
\end{definition}

Figure \ref{show_case} compares the value of $SC_{adj}$ with the original $SC$. The original $SC$ overlooks the topological structure of the Mapper graph, while the $SC_{adj}$ can give a reasonable evaluation. In the next section, we show more examples to validate the effectiveness of $SC_{adj}$. 

\section{Results}\label{Results}
In this section, we compare our proposed D-Mapper algorithm with the classic Mapper algorithm via several experiments. In this work, we implement D-Mapper by expanding the Mapper algorithm in the KeplerMapper package version 2.0.1 \cite{KeplerMapper_JOSS} in Python version 3.11.0. The extended persistence diagram and bottleneck bootstrap are computed by a Python package GUDHI version 3.8.0 \cite{gudhi:urm}.

We compare our proposed D-Mapper algorithm with the classic Mapper algorithm using the metrics $SC_{adj}$ and $SC_{norm}$. We set the interval number $n$ to be identical for both algorithms and use the same clustering algorithm within each hypercube. The clustering method is the density-based spatial clustering of applications with noise (DBSCAN) implemented in the scikit-learn library version 1.1.3. We use a fixed grid to find the best model concerning the $SC_{adj}$ for both the D-Mapper and the classic Mapper. $50$ equally spaced grids in range $(0,\alpha')$ are used to select parameter $\alpha$ in the D-Mapper. Similarly, $50$ equally spaced grids in range $(0,0.5)$ are used to tune the overlapping percentage $p$ in the classic Mapper algorithm. The bottleneck bootstrap sampling steps are set to $100$ to compute the $85\%$ quantile. In this section, we compare the D-Mapper with the classic Mapper via several examples.

\subsection{Two disjoint circles}
\label{Two disjoint circles}
This dataset is created by sampling points randomly from two disjoint circles with centers $(0, 0)$ and $(3, 0)$ and a radius of $1$. $5000$ points are uniformly sampled on each circle. The upper panel of Figure \ref{cirs_results} (a) provides a visualization of the sampled data points. This dataset has a distinctive shape and provides a straightforward performance comparison between the classic Mapper and D-Mapper algorithms.  

\begin{figure}[h]
\centering
    \includegraphics[width=1\textwidth]{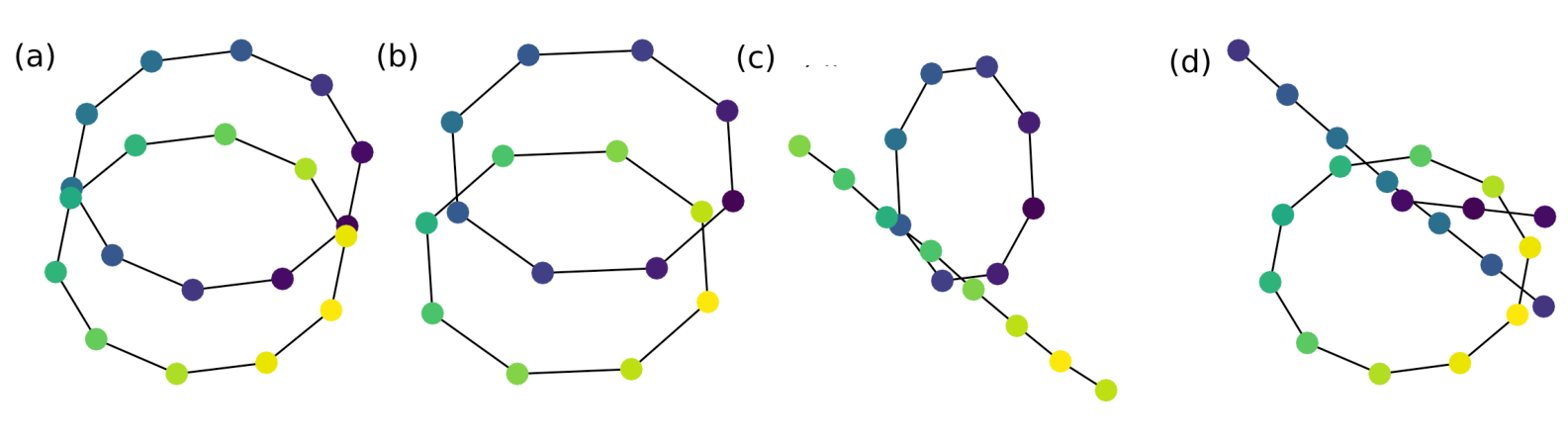}
    \vspace{-2mm}
    \caption{Results of the classic Mapper and D-Mapper on the two disjoint circles. (a) The output graph of the classic Mapper
    with the largest $SC_{adj}$ (\textit{the $1$st row of Table \ref{tab1}}):  $n=12, p = 0.02$.  (b) The output graph of the D-Mapper with the largest $SC_{adj}$  (\textit{the $2$nd row of Table \ref{tab1}}): $n=12, \alpha = 0.127$. (c) An example produced by the classic Mapper with larger $SC$ but lower $TSR$ (\textit{the $3$rd row of Table \ref{tab1}}): $n=12, p = 0.005$.  (d) An example produced by the D-Mapper with larger $SC$ but lower $TSR$ (\textit{the $4$th row of Table \ref{tab1}}): $n=12, \alpha = 0.159$.  }    
    \label{result1}
\end{figure}

In this experiment, we choose the filter function to be a function that projects the original data onto the $X$-axis. We set the number of intervals $n$ to $12$.  The clustering algorithm is DBSCAN with a radius of $0.5$ and a minimum of samples $3$. The parameters $p$ of the classic Mapper algorithm and parameter $\alpha$ of the D-Mapper are tuned via grid search, $p=0.02$ and $\alpha = 0.127$. The comparison of different evaluation metrics is shown in Table \ref{tab1} and Figure \ref{result1}. The resulting intervals are given in  the lower panel of Figure \ref{cirs_results} (a). In Figure \ref{result1}, the color of a node indicates the average value of the projected data within the node. The output graphs of the classic Mapper and D-Mapper are similar (Figure \ref{result1} a-b), and both algorithms capture the topological features of the dataset effectively without noise, resulting in a $TSR$ of 1 for both algorithms. The D-Mapper performs better than the classic Mapper in terms of clustering, as evidenced by its higher $SC_{norm}$. This indicates that D-Mapper has an advantage in identifying more meaningful clusters. Moreover, we display cases in both the classic Mapper and D-Mapper algorithms concerned with the $SC$ only in Figure \ref{result1} (c) and (d). In these two cases, the $SC_{norm}$ values are higher but the $TSR$ values are lower than the results concerned with the $SC_{adj}$. The topology structures of these two cases are different from the dataset in these two figures. The distinct results in this example validate the utility of our proposed metric. It provides a quantitative approach for evaluating Mapper-type algorithms from both topological signal preservation and clustering.

\begin{table}[h]
\caption{Results of the D-Mapper and classic Mapper on the two disjoint circles dataset. The $1$st and $2$nd rows are results with the larger $SC_{adj}$. The $3$rd and $4$th rows are cases when the output graphs have larger $SC$ values but lower $TSR$s.\label{tab1}}%
\begin{tabular*}{\textwidth}{@{\extracolsep\fill}cccl}
\toprule
Algorithm & $SC_{norm}$  & $TSR$ & $SC_{adj}$\\
\midrule
Classic Mapper    &  0.640 & 1.00 & 0.820  \\
D-Mapper    & 0.716 & 1.00 & 0.858  \\
Classic Mapper    & 0.642 & 0.40 & 0.521  \\
D-Mapper    & 0.733 & 0.33 & 0.533  \\
\botrule
\end{tabular*}
\end{table}

\subsection{Two intersecting circles}

\begin{figure}[h]
\centering
    \includegraphics[width=1\textwidth]{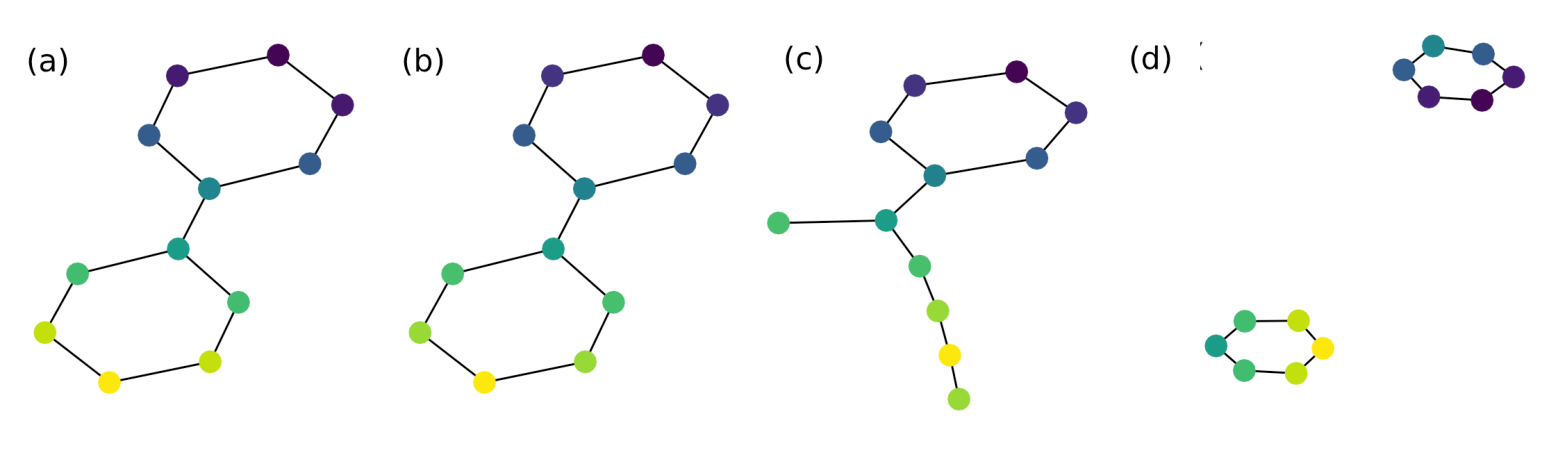}
    \vspace{-2mm}
    \caption{Results of the classic Mapper and D-Mapper on the two intersecting circles. (a) The output graph of the classic Mapper with the largest $SC_{adj}$  (\textit{the $1$st row of Table \ref{tab2}}): $n=8, p = 0.02$. (b) The output graph of the D-Mapper with the largest $SC_{adj}$ (\textit{the $2$nd row of Table \ref{tab2}}):  $n=8, \alpha = 0.088$. (c) An example produced by the classic Mapper with larger $SC$ but lower $TSR$ (\textit{the $3$rd row of Table \ref{tab2}}): $n=8, p = 0.02$. (d) An example produced by the D-Mapper with larger $SC$ but lower $TSR$ (\textit{the $4$th row of Table \ref{tab2}}): $n=8, \alpha = 0.12$.  } 
    \label{result2}
\end{figure}

\begin{figure}[h]%
\centering
\includegraphics[width=0.85\textwidth]{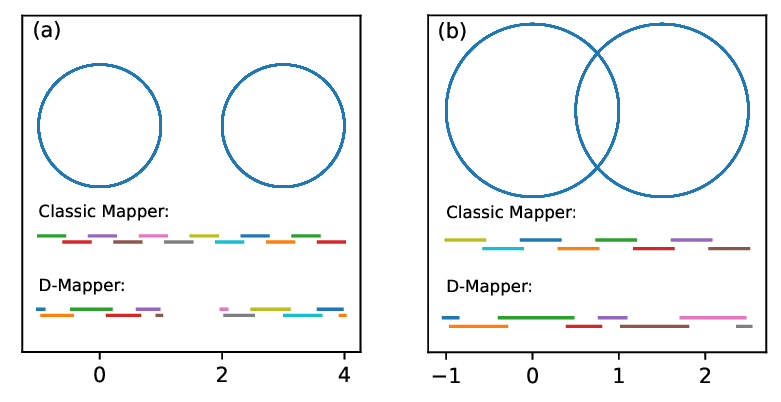}
\vspace{-2mm}
\caption{(a) Upper panel: visualization of two disjointed unit circles whose centers are (0,0) and (3,0), respectively. Bottom panel: the intervals produced by both the classic Mapper and D-Mapper. The intervals produced by the D-Mapper have a gap between two circles, these intervals are more reasonable than  the classic Mapper's. As there are no points between the two disjointed circles, the intervals between the two circles are meaningless. (b) Upper panel: a visualization of two intersecting unit circles dataset whose centers are $(0,0)$ and $(1.5,0)$, respectively. Bottom panel: the intervals produced by both the classic Mapper and D-Mapper. The intervals of the D-Mapper can assign more reasonable intervals than the classic Mapper on the intersecting part.}
\label{cirs_results}
\end{figure}

In this section, we compare the performance of the D-Mapper and the classic Mapper on a two intersecting circles dataset shown in the upper panel of Figure \ref{cirs_results} (b). The data points are generated from two intersecting circles with a radius of $1$ and centers $(0,0)$ and $(1.5,0)$,  respectively. The data generating process is similar to the previous example, $5000$ points are sampled from each circle. The filter function  and clustering algorithm are the same as in the previous example. The number of intervals $n$ is set to $8$. The parameters $p$ of the classic Mapper algorithm and parameter $\alpha$ of the D-Mapper are tuned to $p=0.02$ and $\alpha = 0.088$. The results are given in Figure \ref{result2} and Table \ref{tab2}. The resulting intervals are given in  the lower panel of Figure \ref{cirs_results} (b). As shown in Figure \ref{result2} (a,b), the D-Mapper outperforms the classic Mapper concerning the $SC$ scores (or the $SC_{adj}$) and has the same topological signal rate as the classic Mapper. Similar to the two disjoint circles example, Figure \ref{result2} (c) and (d) also give cases when the metric $SC$ fails to work.

\begin{table}[h]
\caption{Results of the D-Mapper and classic Mapper on the two intersecting circles dataset. The $1$st and $2$nd rows are results with the larger $SC_{adj}$. The $3$rd and $4$th rows are cases when the output graphs have larger $SC$ values but lower $TSR$s.\label{tab2}}%
\begin{tabular*}{\textwidth}{@{\extracolsep\fill}cccl}
\toprule
Algorithm & $SC_{norm}$  & $TSR$ & $SC_{adj}$\\
\midrule
    Classic Mapper &  0.574 & 1.00 & 0.787\\
        D-Mapper& 0.640 & 1.00 & 0.820\\
    Classic Mapper & 0.577 & 0.25 & 0.414\\
    D-Mapper& 0.662 & 0.50 & 0.581\\
\botrule
\end{tabular*}
\end{table}

Both the results of the two disjoint circles and the two intersecting circles indicate that our proposed metric $SC_{adj}$ is more stable than the metric $SC$ measuring both the quality of overlap clustering and extended persistent homology of output of Mapper type algorithm. 

\bigskip

\subsection{3D cat dataset}
\begin{figure}[h]
    \centering
    \includegraphics[width=0.8\textwidth]{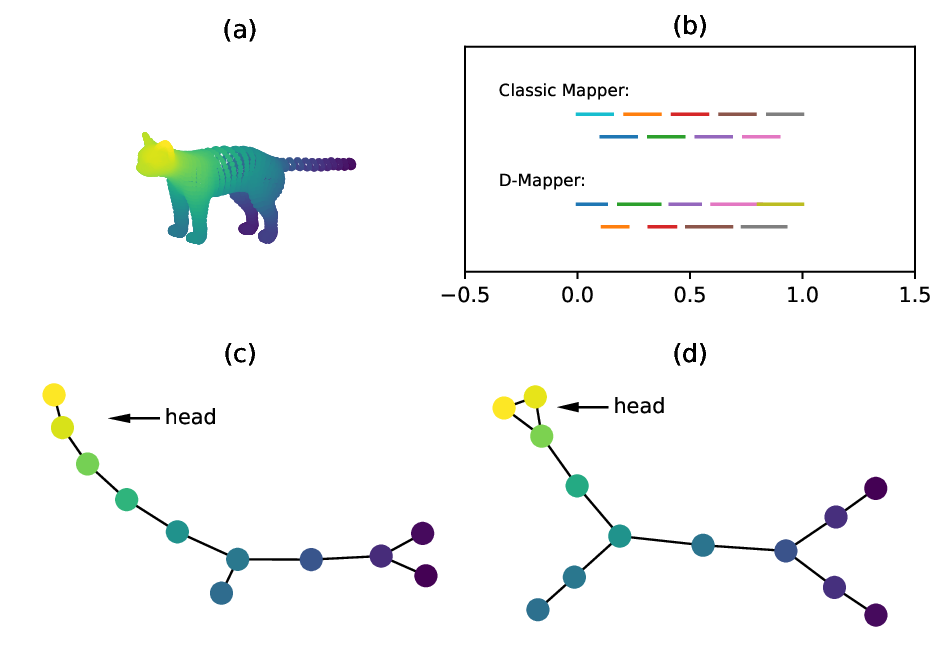}
    \caption{Results of the classic Mapper and D-Mapper on the 3D cat dataset. (a) Visualization of the 3D cat dataset. (b) The intervals produced by both the classic Mapper and D-Mapper. The intervals of both the classic Mapper and D-Mapper are very similar except for intervals near  $1$; this difference gives more details on the cat head in D-Mapper. (c) The graph of the classic Mapper with the largest $SC_{adj} (0.740)$: $n = 9, p = 0.32$. (d) The graph of the D-Mapper with the largest $SC_{adj} (0.755)$: $n=9, \alpha = 0.01$.}
    \label{result3}
\end{figure}

In this experiment, we compare the D-Mapper algorithm with the classic Mapper algorithm on a 3D cat dataset with a more complex topology structure. The 3D cat dataset is originally created by \cite{2004Deformation}. Figure \ref{result3} (a) shows a visualization of the 3D cat dataset. The filter function is the sum of all coordinates. The parameters of the classic Mapper are set to $n=9,~p=0.32$; and the parameters of the D-Mapper are set to $n=9,~\alpha = 0.01$. The clustering method is DBSCAN with a radius of $0.1$ and a minimum of samples $5$. The results and evaluation metrics are shown in Figure \ref{result3} (b,c,d) and Table \ref{tab3}. 

The D-Mapper algorithm results in higher $SC_{norm}$ scores and the same $TSR$ compared to the classic Mapper algorithm. {The output graph of the D-Mapper has one more loop than the classic Mapper, which represents the relationship between the ears, face (including jaw) and neck of the cat (yellow, olive green and light green dots in Figure \ref{result3} c-d).} This example also shows that the D-Mapper could capture subtle features of complex objects better than the classic Mapper due to the flexible intervals and overlapping ratios. 

 \begin{table}[h]
\caption{Results of the classic Mapper and D-Mapper on the 3D cat dataset.\label{tab3}}%
\begin{tabular*}{\textwidth}{@{\extracolsep\fill}cccl}
\toprule
Algorithm & $SC_{norm}$  & $TSR$ & $SC_{adj}$\\
\midrule
    Classic Mapper &  0.480 & 1.0 & 0.740\\
    D-Mapper& 0.510 & 1.0 & 0.755\\
\botrule
\end{tabular*}
\end{table}

\subsection{Covid-19 dataset}
In this section, we apply our proposed method to a real dataset. The Covid-19 pandemic poses an enormous threat to global health and economics. To combat the epidemic, one of the fundamental tasks is identifying and monitoring the mutation of the SARS-COV-2 coronavirus. Studying virus evolution helps understand the basic biological characteristics, develop vaccines and drugs, and forecast trends. A traditional way to study viral mutations is to construct phylogenetic trees based on genetic sequences. However, tree structure representations focus on  vertical processes. Whereas, some processes are horizontal. For viruses,  homologous recombination and reassortment are typical ways that can lead to non-tree-like representation \cite{chanTopologyViralEvolution2013}. We employ the D-Mapper algorithm to analyze a dataset of SARS-CoV-2 coronavirus RNA sequences, utilizing it as a rapid exploratory data analysis tool. This approach allows us to examine viral mutations' vertical and horizontal evolutionary processes across the entire genomic scale. The dataset contains $357$ SARS-COV-2 coronavirus whole RNA sequences with different lineages from \href{https://ngdc.cncb.ac.cn/genbase/}{GenBase} (https://ngdc.cncb.ac.cn/genbase/) in National Genomics Data Center \cite{cncb-genbase2023}. 

We first compute the distance matrix of input sequences. There are several methods available to measure the distance between RNA or DNA sequences, such as alignment algorithm-based \cite{kumar2007multiple} or likelihood-based method \cite{wu2001statistical}, however these methods are computational expensive for large whole DNA or RNA sequences dataset \cite{pearson2013introduction}. To simplify computation, we use $K$-mers \cite{K-mers} frequency vector as a representation of every sequence and compute pairwise distance based on these vectors. We set $k=3$, then every sequence is a $64$ dimension vector. Due to the high similarity of these RNA sequences, a min-max scaling is performed on the distance matrix. The mean value function is selected as the filter function.


The number of intervals are set to  $n = 15, \alpha=0.006$ via grid tuning in the D-Mapper algorithm. The DBSCAN with a radius of $0.6$ and a minimum of samples $3$ is chosen in the clustering algorithm. The results are shown in Figure \ref{covid19}, and the full list is provided in the supplementary file. The output structure represents the evolutionary process of these lineages. The loop structure represents the horizontal evolutionary process. The lineages of the two isolated nodes may indicate significant differences from others,  these lineages may warrant further investigation. In addition, the result of the classic Mapper algorithm is provided in Figure \ref{covid19_C} in Appendix \ref{classic_covid19}.

\begin{figure}[h]%
\centering
\includegraphics[width=0.7\textwidth]{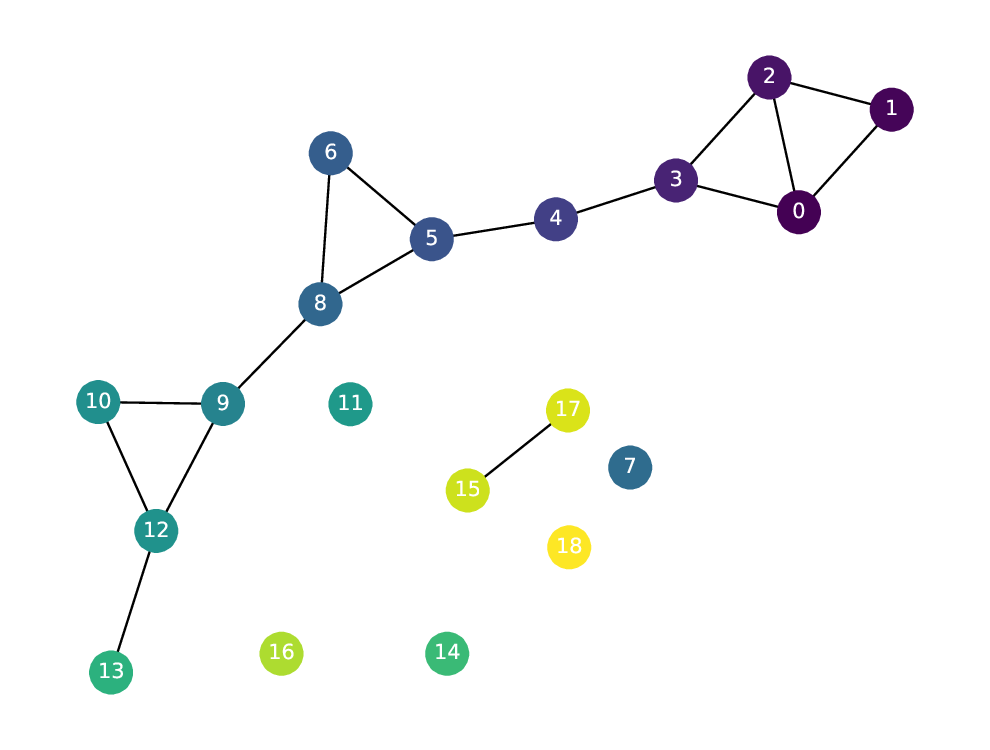}
\vspace{-4mm}

\caption{The result of D-Mapper on the SARS-COV-2 dataset. There are $19$ nodes in total, indexed from $0$ to $18$.}
\label{covid19}
\end{figure}


Some important lineages in these nodes are shown in Table \ref{tab4}. The bottom green loop nodes (nodes 9, 10, 12) contain many XBB variants, which form the recombination of two subvariants and cause large-scale infections worldwide \cite{VirologicalCharacteristicsSARSCoV22023}. {Nodes 4 and 5 contain the JN.1 variant, which is currently widespread worldwide \cite{looi2023covid}. When compared to its predecessor, BA.2.86, the JN.1 variant only has one additional mutation within its spike protein. Nonetheless, BA.2.86 descent while JN.1 rapidly become predominant \cite{wang2024sars}. The D-Mapper graph can give some insights into this phenomenon, where BA.2.86 and JN.1 are both in a loop structure (nodes 5, 6, 8), which indicates the existence of complex horizontal evolution between these two variants. Moreover, JN.1 variant is located in both linear (node 4) and cyclical structures (node 5), suggesting it may play a pivotal role in the evolutionary process. Consequently, other variants presented at these nodes (nodes 4, 5) warrant further scrutiny. The upper purple loop (nodes 0, 1, 2, 3) encompasses variant BA.2, which is the predecessor of BA.2.86. The spike protein of BA.2.86 has more than 30 mutations than BA.2 \cite{planas2024distinct}, which are reflected as significant distances on the D-Mapper graph. In summary, the D-Mapper algorithm offers crucial insights into the rapid and complex evolutionary process of viruses.}  

 \begin{table*}[h]
\caption{Lineages in some important nodes. Only some important shared lineages in loop nodes are shown in this table, the complete information can be found in the supplementary materials.\label{tab4}}%
\begin{tabular*}{\textwidth}{@{\extracolsep{\fill}}lcccccc@{\extracolsep{\fill}}}
\toprule
Nodes & Lineages\\
\midrule
    9, 10, 12 & BA.5.1.3, XBB.1.9.2, XBB.1.42, XBB.1.16.7, XBB.1.17.1\\
    5, 6, 8 & BA.2.86, FY.1.1, XBB.1.12, XBB.1.5.32, FL.2.3.1, JN.1\\
    0, 1, 2, 3& BA.2, EG.4.2, EG.5.2, FY.2, FL.1, XBB.2.3.8 \\
\botrule
\end{tabular*}
\end{table*}

\section{Discussion}\label{Discussion}
Our proposed algorithm is a probabilistic model-based approach that uses data intrinsic characteristics and probability models to generate distribution-guided covers and improved topological features. It is a viable alternative to non-probabilistic approaches. With the D-Mapper, once we get the distribution of the projected data, we can further get the intervals from each distribution. Theoretically, the distribution of original data can be computed by transformation of distributions if the filter function is well proposed, this may be helpful for further theoretical analysis.    

The selection of the optimal number of components remains a difficult problem to address. One possible solution is to use a non-parametric mixture model, such as the Dirichlet process, by introducing the Dirichlet process prior to choose the number of components adaptively. Another approach is to apply the information criteria like AIC or BIC to find a proper number of components for the mixture model.
 
 Although our proposed metric $SC_{adj}$ is more reasonable than $SC_{norm}$, it still has some limitations. The $TSR$ primarily assesses the stability of the Mapper graph but does not account for other aspects of the Mapper graph, for example, the richness of the topological structures. Moreover, domain knowledge is essential for interpreting a Mapper graph. It is advisable to integrate domain expertise when assessing different Mapper graphs in real applications. 

\section{Conclusion}\label{Conclusion}

In this work, we propose a distribution-guided Mapper algorithm (D-Mapper) to relax the fixed intervals and overlapping ratios restriction in the classic Mapper algorithm.  Our proposed algorithm combines a mixture model and the Mapper algorithm to obtain irregular intervals based on the density of the projected data. With these irregular intervals, the D-Mapper algorithm can gain deeper topological insights and enhance clustering outcomes. To validate the effectiveness of our proposed algorithm, we introduce the $SC_{adj}$ score to combine $SC$ scores and extended persistence diagram as a metric to reflect the performance of overlap  clustering and persistent homology. 

We also conduct simulation studies with different complexity to evaluate the D-Mapper algorithm by the $SC_{adj}$. In all simulations, the D-Mapper outperforms the classic Mapper algorithm concerning the metric $SC_{adj}$. The $SC_{norm}$ of the D-Mapper are all higher than the classic Mapper algorithm, and all $TSR$s are $1$. This indicates that the D-Mapper can achieve better clustering than the classic Mapper algorithm while outputting high-quality Reeb graph approximations. Note that we tune parameters by grid searching the highest $SC_{adj}$, thus the $TSR$ is relatively high. In many cases, the $TSR$ could be low, as we show in Table \ref{tab1} and Table \ref{tab2}. Our experimental results also indicate that our proposed metric $SC_{adj}$ is more stable than the metric $SC$. We also apply the D-Mapper algorithm to the SARS-COV-2 coronavirus RNA sequences, and the result shows that the D-Mapper algorithm can reflect both the vertical and horizontal evolutionary processes.

\begin{appendices}

\section{The bottleneck bootstrap algorithm} \label{bootstrap}
The confidence sets for extended persistence diagrams allow us to separate topological signals from noises. One effective method to compute confidence sets is the bottleneck bootstrap algorithm \cite{fasyConfidenceSetsPersistence2014}. The bottleneck distance is often used as a metric to compute the distance between persistence diagrams. For an extended persistence diagram, points that have a short lifetime are usually considered as a noise and a long lifetime could be a signal. 

Paired the original data points and the filtered values as $(X_1, Y_1),...,(X_n, Y_n)$. The Mapper graph is represented by $M$ and the bootstrap samples are denoted as $(X^*_1, Y^*_1),...,(X^*_n, Y^*_n)$. The Mapper graph generated from these bootstrap samples is denoted by $M^*$. Given a confidence level $\epsilon$, this algorithm outputs a distance $d_{\epsilon}$, which can be utilized to differentiate between noise and signal in the extended persistence diagram. Algorithm \ref{confidence_sets} gives the details of the bottleneck bootstrap algorithm.

\begin{algorithm}[H]
    \caption{The bottleneck bootstrap algorihtm}
    \label{confidence_sets}
    \begin{algorithmic}[1]
    \State Input: Paired orignial data and filtered values $(X_1, Y_1),...,(X_n, Y_n)$, all parameters for the Mapper or D-Mapper algorithm, a confidence level $\epsilon$.
    \State Generating a Mapper graph $M$ by the Mapper or D-Mapper algorithm by using original samples.
    \For{$i=1$ to N}
    \State Bootstrap sampling: $(X^*_1, Y^*_1),...,(X^*_n, Y^*_n)$.
    \State Generating a bootstrap Mapper graph $M^*$ by the Mapper or D-Mapper algorithm.
    \State $d_i = d(M,M^*)$
    \EndFor
    \State Find the $\epsilon$ quantile from $\{d_1,...,d_N \}$, denote as $d_{\epsilon}$.
    \State Return $d_{\epsilon}$.
    \end{algorithmic}
\end{algorithm}

\section{The output of the classic Mapper algorithm on the Covid-19 dataset}\label{classic_covid19}

The result of the classic Mapper algorithm is shown in Figure \ref{covid19_C}, and the full list is provided in the supplementary file.  

\begin{figure}[H]%
\centering
\includegraphics[width=0.7\textwidth]{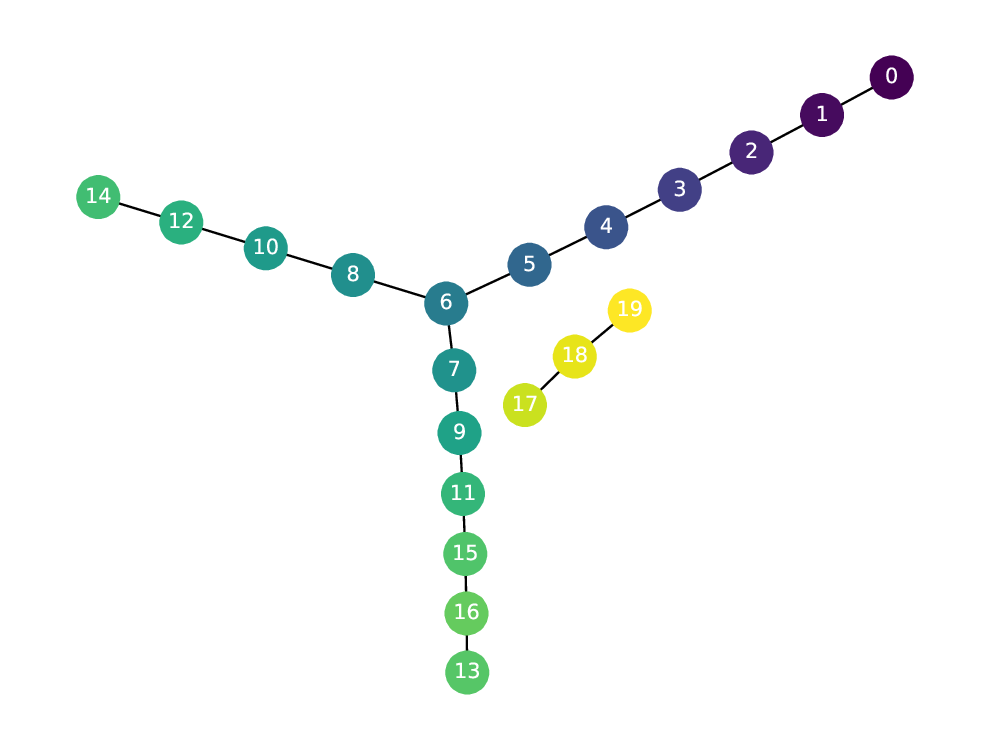}
\vspace{-4mm}
\caption{The result of the classic Mapper on the SARS-COV-2 dataset  ($n=15, p=0.49$). There are $20$ nodes in total, indexed from $0$ to $19$.}
\label{covid19_C}
\end{figure}






\end{appendices}

\subsection*{Acknowledgements}

The authors thank ShanghaiTech University for supporting this work through the startup fund and the HPC Platform. 

\subsection*{Authors' contributions}
YT—conceptualization, data curation, formal analysis, methodology, coding, visualization, writing—original draft. SG—conceptualization, data curation, formal analysis, methodology, coding evaluation, project administration, supervision, writing—final draft, writing—review and editing.  All authors read and approved the final manuscript.

\subsection*{Funding}
YT and SG are supported by the Shanghai Science and Technology Program (No. 21010502500), the National Natural Science Foundation of China (12401383).

\subsection*{Availability of data and materials}

The paper employs publicly accessible datasets in the 3D cat and Covid-19 experiments. The simulated data of the two disjoint circles and two intersecting circles, the results of the all experiments and code are available in the public repository: https://github.com/ShufeiGe/D-Mapper. 

\section*{Declarations}

\subsection*{Ethical approval and consent to participate}
Not applicable.

\subsection*{Consent for publication}
Not applicable.

\subsection*{Conflict of Interests}
The authors declare that they have no competing interests.



\bibliography{sn-bibliography}

\end{document}